\documentclass[a4paper,12pt]{article}
\usepackage[latin1]{inputenc}
\usepackage{amsmath,amssymb,amsthm,cite,verbatim}
\usepackage[all]{xy}
\usepackage[bookmarks=true]{hyperref}
\usepackage{color}
\usepackage{graphicx}
\usepackage{hyperref}
\title{Bielliptic modular curves $X_0^*(N)$ with square-free levels}

\author{ \ Francesc Bars\footnote{First author is supported by MTM2016-75980-P and MDM-2014-0445} and Josep Gonz\'alez  \footnote{The  second author is partially supported by DGI grant  MTM2012-34611.}}

\newcounter{ap}
\newtheorem{prop}[ap]{Proposition}
\newtheorem{lema}[ap]{Lemma}
\newtheorem{teo}[ap]{Theorem}
\newtheorem{cor}[ap]{Corollary}
\newtheorem{rem}[ap]{Remark}

\theoremstyle{definition}

\theoremstyle{remark}

\numberwithin{equation}{section}

\newcommand{\Q}{\mathbb{Q}}

\newcommand{\Z}{\mathbb{Z}}

\newcommand{\F}{\mathbb{F}}

\newcommand{\PP}{\mathbf{P}}
\newcommand{\C}{\mathbb{C}}

\newcommand{\Gal}{\mathrm{Gal}}

\newcommand{\End}{\operatorname{End}}
\newcommand{\Aut}{\operatorname{Aut}}

\newcommand{\Frob}{\operatorname{Frob}}
\newcommand{\Jac}{\operatorname{Jac}}

\newcommand{\New}{\operatorname{New}}

\newcommand{\cI}{{\mathcal I}}

\newcommand{\cL}{{\mathcal L}}

\date{}
\begin{document}
\maketitle

\begin{abstract}
\noindent Let $N\geq 1$ be a square-free integer such that the
modular curve $X_0^*(N)$ has genus $\geq 2$. We prove that
$X_0^*(N)$ is bielliptic exactly for  $19$ values of $N$, and we
determine the automorphism group of these bielliptic curves. In
particular, we obtain the {{first}} examples of nontrivial
$\Aut(X_0^*(N))$ when the genus of $X_0^*(N)$ is $\geq 3$. Moreover,
we prove that the set of all quadratic points over $\mathbb{Q}$ for
the modular curve $X_0^*(N)$ {{with genus $\geq 2$ and $N$
square-free}} is not finite exactly for $51$ values of $N$.

\end{abstract}

\section{Introduction}
Let  $X$ be a smooth projective curve defined over a number field
$K$ of genus $g_X$ at least
 two. In \cite{Fa},   Faltings proved the finiteness of the set of  points
 of $X$ defined over $K$, denoted by $X(K)$.
 After that, for a finite extension $L/K$, the natural object to consider was the set
of  points of $X$ defined over all quadratic extensions of  $L$,
i.e. the set
$$\Gamma_2(X,L):=\cup_{[F:L]\leq 2} X(F).$$
 In \cite{SiHa}, Harris and Silverman proved that the above set is not
finite for some number field $L$ if, and only if, $X$ is
hyperelliptic or bielliptic, i.e.  the curve
$X\times_K\overline{K}$ admits a degree 2 map to the projective line
or to  an elliptic curve over
 a fixed algebraic closure $\overline{K}$ of $K$.
Moreover,  from the work of Abramovich, Harris and Silverman,  in
\cite[Theorem 2.14]{BaMom} it is proved that the set
$\Gamma_2(X,L)$ is infinite if, and only if,  $X$ is hyperelliptic
over $L$, i.e. there is a morphism $ \phi\colon
X\rightarrow\mathbb{P}^1_L$ of degree two defined
 over $L$, or $X$ is bielliptic over $L$, i.e.
there exist an elliptic curve $E$ over $L$ and a morphism
$\phi:X\rightarrow E$ of degree two defined over $L$, such that the
$L$-rank of $E$ is at least one.

This made   the study of  bielliptic curves a matter of deep interest for
Arithmetic Geometry. This was developed for the modular world,
 {{because}}  $L$-points of modular  curves have a moduli interpretation on
elliptic curves. The first work concerned the modular curves
$X_0(N)$. The levels $N$ for which
   the set  $\Gamma_2(X_0(N),\mathbb{Q})$ is finite are determined in \cite{Ba1}.

 Later, different results determining the bielliptic curves among  some
modular curves recovering $X_0(N)$ can be found in \cite{JK1} for
$X_1(N)$, in \cite{JKS} for $X_{\Delta}(N)$  and in \cite{BKX} and
\cite{JK2} for $X(N)$.

Two important tools to obtain such results are the following. First,
if there is a  morphism of curves $X\twoheadrightarrow X'$ such that
$X$  is bielliptic  and the genus of $X'$ is at least $2$, then $X'$
is hyperelliptic or  bielliptic \cite[Proposition 1]{SiHa}. One can
use results about hyperelliptic modular curves, whose study has been
widely treated in the last decades. Second, $X$ is bielliptic if,
and only if, there exists an involution of $X$ fixing
$2g_X-2$-{{many}} points, where $g_X$ denotes de genus of $X$.

In this papr we consider the modular curves  $X_0^*(N)$. They are
defined as the quotient of the modular curve  $X_0(N)$ by the  group
of  all Atkin-Lehner involutions, which is  defined  over $\Q$. Its
$K$-points, which are not cuspidal, correspond to $K$-curves with
additional data from the level $N$. See \cite{Elkies} for further
information.

Here, we restrict our attention to the case where  $N$ is
square-free. Under such assumption, the modular curve $X_0^*(N)$
corresponds to the quotient $X_0(N)/\Aut(X_0(N))$  and, thus,  it
does not have any natural automorphisms  (in particular,
involutions) coming from $X_0(N)$, except for  $N=37$. { { These
curves have two properties that play an important role in the
development of this article. On the one hand, all involutions of
$X_0^*(N)$ are defined over $\Q$. On the other hand, the
endomorphism algebra $\End(\Jac(X_0^*(N))\otimes \Q$ is isomorphic
to the product of totally real numbers fields (cf. \cite[\S 2]{BH}).
Any of these properties can fail when $N$ is non square-free and the
study of this case needs additional  tools.}}

We point out that if $X_0^*(N)$ is bielliptic, then {{ a bielliptic
quotient
 (i.e. the quotient  of $X_0^*(N)$ by  a bielliptic
involution),}} is an elliptic curve $E$ defined over $\Q$ of
conductor $M|N$ with odd analytic rank, because the attached newform
is invariant under the Atkin-Lehner involution $w_M$. Hence, it is
expected that the algebraic rank of $E$ is odd.


In our case, the knowledge of the values of $N$
for which $X_0(N)$ is bielliptic is not useful to obtain  bielliptic
curves $X_0^*(N)$, because, in this case, these curves  have genus
at most one or are hyperelliptic.
A new approach is needed to deal with our case and, here,
we use a  method given in \cite{Go17} to discard
automorphisms of certain order
 in the automorphism group of a curve defined over a finite field. In particular, this  method allows us to deduce
 that the automorphism group is trivial for such values of $N$, when such method works.
This approach behaves   well  for  odd square-free integers $N$,
when  $N$ is the product of two {{or}} three primes.
 When it fails, we use the usual method of
reducing modulo a prime $p$ to discard some situations. For the
remaining cases, using {{a}} Theorem of Petri {{\cite{Pe1923}}}, we
implement a method to recognize whether a curve $X_0^*(N)$ is
bieliptic and compute equations for the elliptic quotient.

The main result of this article is the following.
\begin{teo}\label{main} Let $N>1$ be a square-free integer. Assume that the genus of the modular curve $X_0^*(N)$
is at least $2$. Then, the modular curve $X_0^*(N)$ is bielliptic
if, and only if, {{$N$ is in the following table
$$ \begin{array}{c|r|}
\text{genus} & N\phantom{cccccccccccccccccccccccc}\\\hline
2& 106,122,129, 158,166,215,390\\ \hline
3&  178,183, 246, 249, 258, 290, 303,318, 430,455,510\\\hline
4 & 370\\ \hline
\end{array}
$$
}}
 For these values of $N$, the automorphism group of
$X_0^*(N)$  has order $2$ when its genus is greater than two,
otherwise it is the Klein group.
\end{teo}
 Concerning $\Aut(X_0^*(N))$, with $N$
square-free and  genus $\geq 2$, it is known that it is
an abelian 2-group (cf. \cite{BH}).  Moreover, when $N$
is  prime  this group  is nontrivial
 if, and only if,  the genus of the curve is 2 and, in this case, the  group has order $2$ (cf. \cite[Theorem 1.1]{BH}).
 In fact, it is expected that this group is trivial for almost
 {{all}}
 square-free $N$.

In this paper, we can observe that the Klein group appears naturally
for genus two curves which are also bielliptic (cf. Remark
\ref{rem3}). Moreover, we point out  that a bielliptic curve could
have several involutions and also more than one bielliptic
involution when its genus is $\leq 5$  (cf.
\cite[Prop.2.10]{BaMom}). Nevertheless, this does not happen in our
case {{when the genus is $>2$.}}

As a  by-product  of this work, we  also obtain the following result.
\begin{prop}
 The automorphism group of $X_0^*(N)$ is trivial for the following values of
 $N$:
$$
\begin{array}{l}
185, 202, 259, 262, 267, 282, 301, 305, 310, 354, 393, 394, 395, 399, 426, 427,445,\\
458, 462, 546, 570, 581, {{582}}, 602, 710, 786, 795, 903, 1001,
1015\,.
\end{array}
$$
Moreover, the group $\Aut(X_0^*({{366}}))$ has order $2$, and the
quotient curve has genus $2$.
\end{prop}
For many of these values (see Propositions \ref{prop2} and \ref{prop8}), this result is obtained by using the method
to  discard the existence of involutions, which was mentioned above. A Magma code to be applied in our case can be found in

 \url{http://mat.uab.cat/~francesc/programmesXoestrellaMagma.pdf}

\noindent(this html page, also contains  different codes in Magma for computing
the genus of $X_0^*(N)$ and its $\F_{p^k}$-points). For the remaining values of the above proposition (see Propositions \ref{prop4}, \ref{prop5}, \ref{prop10} and  \ref{prop11}),   Theorem of Petri  is the main tool.

As for  quadratic points,  by the work of Hasegawa and Hashimoto
(cf. \cite{HaHa}),  we know that $X_0^*(N)$ is hyperelliptic with
$N$ square-free if, and only if, the curve  has genus 2.  When $X_0^*(N)$ is bielliptic with genus $>2$, the rank of the elliptic quotient  turns out  to be one. So  we conclude
\begin{teo}\label{main2}Let $N>1$ be a square-free integer. Assume that the genus of the modular curve $X_0^*(N)$
is $\geq 2$. Then, the set  $\Gamma_2(X_0^*(N),\mathbb{Q})$ is infinite
if, and only if, $N$ lies in the set
$$\begin{array}{l}
\{{{67}}, 73,85,93,103,106,107,115,122,129,133,134,146,154,158, 161,
165, 166, 167,
 170, 177,178,\\183,
186, 191,205, 206, 209, 213,215, 221, 230, 246, 249, 255, 258, 266,
285, 286, 287, 290, 299, 303,\\
318, 330, 357, 370, 390, 430, 455, 510\}\,.
\end{array}$$
\end{teo}
{{Theorem \ref{main2} holds when we replace
$\Gamma_2(X_0^*(N),\mathbb{Q})$ with $\Gamma_2(X_0^*(N),K)$, where
$K$ is any number field. This is due to the fact that if $X_0^*(N)$
is hyperelliptic or bielliptic over $K$, then it is hyperelliptic or
biellitic over $\Q$ (cf. Lemma \ref{lem1})}}


\section{Preliminary results}
Let $N>1$ be an integer. We fix once and for all the following
notation. We denote by
 $g_N$ and $g_N^*$  the genus of $X_0(N)$ and $X_0^*(N)$, respectively, and $n$ is the number of primes dividing $N$.
 {{For any $1\leq d|N$ with $(d,N/d)=1$ we have an
 involution $w_d\in \Aut(X_0(N))$, called the Atkin-Lehner involution attached to $d$,
 and we denote by $B(N)$ the group of all Atkin-Lehner
 involutions.}}
 We denote by $\New_N$ the set of normalized newforms in $S_2(\Gamma_0(N))$,
 and $\New_N^*$ is the subset of $\New_N$ consisting of  the newforms  invariant under the action of the group  $B(N)$.
 For an integer  $m\geq 1$ and a newform $f\in\New_N$, $a_m(f)$  is the $m$-th Fourier coefficient of $f$. {{For an eigenform $g\in S_2(\Gamma_0(N))$,  $A_g$ denotes  the abelian variety defined over $\Q$ attached by Shimura to $g$. As usual,  $\psi$
 is the Dedekind psi function.}}

In the sequel,  $N$ is  square-free. We recall the following result
of  Baker and Hasegawa.
\begin{lema}[Corollary 2.6 in \cite{BH}]\label{lem1}
 The group $\Aut X_0^*(N)$ is  elementary $2$-abelian and every automorphism of $X_0^*(N)$ is defined over $\Q$.
\end{lema}

From now on, we assume that $X_0^*(N)$ has a bielliptic involution
$u$. Let us denote by $E$ the elliptic quotient $X_0^*(N)/\langle
u\rangle$ and by $\pi$ the nonconstant morphism $X_0(N)\rightarrow
X_0^*(N) \rightarrow E$, which has  degree $2^{n+1}$ and is defined over
$\Q$. Let $M$ be the conductor of $E$. It is well-known that $M|N$
and there exist a morphism $\pi_M\colon X_0^*(M)\twoheadrightarrow E$ and a normalized newform $f_E\in \New_M^*$ such that $\pi_M^*(\Omega^1_{E/\Q})=\Q( f_E(q)
{\rm d\,}q/q)$. Moreover,  $\pi^*(\Omega^1_{E/\Q})=\Q( g(q){\rm d\,}q/q)$, where $g$ is the eigenform
$\sum_{1\leq d| N/M}w_d^*(f)\in S_2(\Gamma_0(N))^{B(N)}$ and $w_d$ stands for the Atkin-Lehner involution attached to $d$. More precisely, $g(q)=\sum_{1\leq d| N/M} d \,f_E(q^d)$.

Since the gonality of $X_0^*(N)$ is $\leq 4$, by applying Proposition 4.4 in \cite{BGGP}, we obtain  $g_N^*\leq 35$. Nevertheless,
this fact is not very useful in order to determine a finite set of possible values for  $N$.
The following lemma helps us to achieve this goal. Note that for a prime $p\nmid M$,
 $|E(\F_{p^2})|\leq (p+1)^2$,  but $|E(\F_{p^2})|= (p+1)^2-a_p^2(f_E)$   due to the congruence of  Eichler-Shimura.

\begin{lema}\label{psi}  The following inequalities  hold:
\begin{itemize}
\item[{\rm (i)}] If $p\nmid N$, then
$$
{\rm (a)}\,\, \displaystyle{\frac{\psi (N)}{2^n}\leq 12\cdot \frac{2|E(\F_{p^2})|-1}{p-1}}\,, \quad
{\rm (b)}\,\, \displaystyle{g_N^* \leq 2 \frac{|E(\F_{p^2})|}{p-1}} \,,\quad
{\rm (c)}\,\,  \displaystyle{g_N \leq  2^{n+1}\frac{|E(\F_{p^2})|}{p-1}}\,.$$

\item [{\rm (ii)}] If     $p\mid N$,  then   $g_{N/p}^*  \leq 1$ or
   $$   \begin{array}{cccc}
\displaystyle{  \frac{\psi (N/p)}{2^{n-1}}\leq   12 \frac{2\,|E(\F_{p^2})|-1}{p-1}}\,,
& g_{N/p} ^*\leq \displaystyle{2\frac{ |E(\F_{p^2})|}{p-1}}\,,& \displaystyle{g_{N/p} \leq  2^{n}\frac{|E(\F_{p^2})|}{p-1}\,, }& \text{if $p\nmid M$,}\\  [ 8 pt]
\displaystyle{ \frac{\psi (N/p)}{2^{n-1}}}\leq   12\frac{2 p^2+1}{p-1}\,, & \displaystyle{g_{N/p}^* \leq 2\frac{p^2+1}{p-1}}\,,
& \displaystyle{g_{N/p}\leq 2^n\frac{p^2+1}{p-1}}\,,&\text{if $p|M$.}\end{array}
$$
 \end{itemize}
\end{lema}

\noindent{\bf Proof.} Assume $p\nmid N$. We generalize the argument
used by Ogg in \cite{Ogg}. Indeed,  $X_0(N)(\F_{p^2})$ contains
$2^n$ cusps and at least $(p-1)\frac{\psi(N)}{12}$ {{many}}
supersingular points (cf.    \cite[Lemma 3.20 and 3.21]{BGGP}).
Since there is a nonconstant morphism defined over $\Q$ from
$X_0(N)$ to an elliptic quotient $E$ of $X_0^*(N)$ which has degree
$2^{n+1}$, $|X_0(N)(\F_{p^2})|\leq 2^{n+1} |E(\F_{p^2})|$. Parts (b)
and (c) in (i) are  obtained applying \cite[Lemma 3.25]{BGGP}.

If $p|N$, then $X_0(N)/\F_p$ is the copy of two curves
$X_0(N/p)/\F_p$ , and  the normalization of $X_0^*(N)/\F_p$ is the
curve  $X_0^*(N/p)/\F_p$ (cf. \cite{FM}). If the reduction of the
involution $u$
 is the identity, then $g_{N/p}^*$ is the genus of $E/\F_p$. Otherwise, $|X_0^*(N/p)(\F_{p^2})|$     is at most $2(p^2+1)$ or
$|E(\F_{p^2})|$,      depending on whether $p\mid M$ or not.
\hfill $\Box$

\begin{rem}
The above conditions imply $n\leq 4$, when $N$ is odd, and $n\leq
5$ in the even case.   The values    $N$ for which $g_N^*\leq 1$ can
be found  in  \cite[Proposition 3.1 and 3.2]{GL},    and those for which $g_N^*=2$
can be found in    \cite[Theorem 2 ]{Hata}.
\end{rem}

Keeping the above notation, we present the following lemma, which will used to discard some  elliptic curves $E$ for a value $N$.
\begin{lema}\label{degree}  Let $E'$ be the elliptic curve in the $\Q$-isogeny class of $E$ that is an optimal quotient of the jacobian of $X_0^*(M)$.
If $M=N$, then the degree $D$  of the modular parametrization $\pi_N \colon X_0(N)\rightarrow E'$ divides $2^{n+1}$.
\end{lema}

\noindent{\bf Proof.} The statement follows from the optimality of $\pi_N$ and the fact that the degree of $\pi$ is $2^{n+1}$.
 \hfill $\Box$

\begin{rem}
The degree $D$  can be found in \cite[Table 5]{Cre}.
\end{rem}

\section{Odd case}  If $N$ is odd, applying   Lemma \ref{psi} for $p=2$, we have $\psi (N)/2^n\leq 204$.
This fact  implies  $n\leq 3$, except for  the  values $N\in\{
3\cdot 5\cdot 7\cdot 11, 3\cdot    5\cdot 7\cdot 13\}$. The case
$n=1$ can be discarded,   since $\Aut (X_0^*(p))$ is trivial for all
$p$ except for {{$g_p^*=2$}} and, in this particular hyperelliptic
case, the automorphism group  has order two (cf. \cite[Theorem
1.1]{BH}). {{We assume also that $g_N^*\geq 2$ and apply Lemma
\ref{psi} for the values $\gcd (N,3)=1 $ for which $\psi(N)/2^n\leq
186$.}}
 There are exactly {{$146$}}
values for odd $N$ such that $n>1$, $\psi(N) \leq 204\cdot 2^n$
{{(or $\leq 186\cdot 2^{n}$ if $6\nmid N$ )}} and all  these cases
satisfy $1<g _N\leq 9\cdot 2^{n+1}$. More precisely, we have
{{$100$}} cases for $n=2$, {{$44$}} for $n=3$ and $2$ for $n=4$.

 We can reduce this list by
  considering the pairs $(N,E)$, where  $E$ is the $\Q$-isogeny class of the elliptic curves  of conductor $M|N$ such that its attached newform
 $f_E$ lies in $\New _M^*$.
 From \cite[Table 5,Table 3]{Cre}, we obtain the degree $D$ of $\pi_N$, when $M=N$, and  $a_2(f_E)$, $a_3(f_E)$ {{and $a_5(f_E)$}}.
 In particular, we know  $|E(\F_4)|$,
$|E (\F_9)|$ when $3\nmid M$, {{and $|E(\F_{25})|$ when $5\nmid
M$}}. For $ n> 1 $, we can discard the pairs $(N,E)$  that do not
satisfy the conditions in Lemmas \ref{psi} and \ref{degree}. { {When
$g_N^*=2$, if an elliptic quotient  of $X_0^*(N)$  is not
bielliptic, then we can discard $N$ (cf. Remark \ref{rem3}). In
particular, we  discard $N=285$ because
 the elliptic quotient of $X_0^*(285)$ with  conductor $285$ does not satisfy Lemma \ref{degree}. }} In Table 1, we
present the remaining possibilities, where the label of the elliptic
curve $E$ is the one in Cremona tables.

$$           \begin{array}{|l|r|r|r|c|}
          \hline N&g_N&g_N^*&M&\text{Label $E$}\\
          \hline 129=3\cdot 43& 13& 2&N&a\\
          &&&43&a\\
          \hline 183=3\cdot 61& 19& 3&61&a\\
          \hline 185=5\cdot 37& 17&  3&37& a\\
          \hline 215=5\cdot 43& 21& 2&N&a\\
                              &   &  & 43 &a\\
          \hline 237=3\cdot 79& 25& 5&79&a\\
          \hline 249=3\cdot 83& 27& 3&N&b\\
          &&&83&a\\
          \hline 259= 7\cdot 37& 23& 4&37&a\\
          \hline 267= 3\cdot 89& 29& 4&89&a\\
          \hline 301=7\cdot 43& 27& 6&43&a\\
          \hline 303=3\cdot 101& 33& 3&101&a\\
          \hline 305=5\cdot 61& 29& 4&61&a\\
                    \hline 393=3\cdot 131& 43& 5&131&a\\
              \hline 395=5\cdot 79& 39& 4&79 &a     \\
              \hline 415=5\cdot 83& 41& 8&83&a\\
              \hline 427=7\cdot 61& 39& 8&61&a \\
              \hline 445=5\cdot 89& 43& 7&89 &a\\
          \hline 581= 7\cdot 83& 55& 8&83&a\\\hline
              \end{array}
 \quad
         \begin{array}{|l|r|r|r|c|}
          \hline
          N & g_N& g_N^*&M& \text{Label $E$}\\ \hline
                                                    \hline 273=3\cdot 7\cdot 13 &33 &4& 91 &a\\
                    \hline 385= 5\cdot 7\cdot 11& 45 & 4 &77& a\\
                  \hline 399=3\cdot 7\cdot 19& 49& 4&57&a\\
         \hline 429=3\cdot 11\cdot 13& 53& 3&143&a\\
         \hline 435=3\cdot 5\cdot 29& 57& 5&145&a\\
         \hline 455=5\cdot 7\cdot 13& 53& 3&65&a\\
                                     &   &   &91 &a\\
         \hline 465=3\cdot 5\cdot 31& 61& 5&155&c\\
                  \hline 555=3\cdot 5\cdot 37 &73& 5& 185&c\\
         \hline 615=3\cdot 5\cdot 41& 81& 6&123&b\\
         \hline 645=3\cdot 5\cdot 43& 85& 5&215&a\\
                                     &  &  &129 & a\\
                                     \hline 705=3\cdot 5\cdot 47 &88&8&141 &d\\
         \hline 715=5\cdot 11\cdot 13 & 81 & 8 &143 &a\\
                                      &    &    & 65 &a\\
         \hline 795=3\cdot5\cdot 53&105&10& 265 &a\\
                &&& 53&a\\
         \hline {{861=3\cdot7\cdot41}}&109&7&123&b\\
         \hline 903=3\cdot 7\cdot 43& 113& 11&129&a\\
         \hline 987=3\cdot 7\cdot 47& 125& 9&141&d\\
         \hline 1001=7\cdot 11\cdot 13& 109& 8&143&a\\
                              \hline 1015=5\cdot 7\cdot 29& 117& 11&145&a\\
                   \hline
             \hline
         1155=  3\cdot 5\cdot 7\cdot 11  &185& 8  & 77& a  \\ \hline
         1365  =3\cdot 5\cdot 7\cdot 13 &237& 9 & 455& a  \\
             & &  & 65 &a \\
 \hline
             \end{array}
             $$
     $$ \text{ Table 1} $$

  \vskip 0.2 cm First, we examine the hyperelliptic cases in Table 1, which correspond to those such  that $g_N^*=2$ (cf. \cite[ Theorem 2]{FM}).

\begin{prop}
The curves of genus two $X_0^*(129)$ and  $X_0^*(215)$   are bielliptic.
\end{prop}
\noindent{\bf Proof.} For $N=129$ and  $N=215$, the jacobian of $X_0^*(N)$ is isogenous over $\Q$ to the product of two elliptic curves
$E_1\times E_2$, where $E_1$ has conductor $N$ and $E_2$ has conductor $M=43$.
 Hence, there
exist two normalized newforms $f_1\in \New_N^*$ and $f_2\in \New_M^*$ such that the elliptic curves $A_{f_1}$ and $A_{f_2}$ are isogenous over
$\Q$ to $E_1$ and $E_2$ respectively.
The set of the   regular  differentials
$$
\omega_1=f_1(q) {\rm d\,} q/q\,,\quad \omega_2=(f_2(q)+ \ell f_2(q^\ell) ){\rm d\,} q/q\,,\text{ with }\ell=N/M\,,
$$
is a basis of $\Omega^1_{X_0^*(N)/\Q}$.
The functions $\displaystyle{x=\frac{\omega_2}{\omega_1}}$ and $ \displaystyle{y= \frac{ {\rm d}\, x}{\omega_1}}$ on $X_0^*(N)$ satisfy the
equations
$$\begin{array}{c |rcr|}
N & & & \text{equations}\qquad\qquad\qquad\qquad\qquad\qquad\\ \hline
129 & 4\, y^2 &= &  x^6-11\,  x^4+35\, x^2-9\\ [6 pt]
215 & 4\,y^2 &= & - x^6-5\, x^4-3\, x^2+25 \\ \hline
\end{array}
$$
For $N=129$ and $ 215$, it is clear that  the curves have two
bielliptic involutions $(x,y)\mapsto (-x,\pm y)$.
\hfill $\Box$

\begin{rem}\label{rem3}
Assume that  $X_0^*(N)$ has genus two {{ and has an elliptic
quotient. Then, $\Jac(X_0^*(N))$ is isogenous over $\Q$ to the
product of two non isogenous elliptic curves $E_1$ and $E_2$. }} If
$X_0^*(N)$ has a bielliptic involution  $u$,  then $u^*(E_i)=E_i$
and their  regular differentials $\omega_1$ and $\omega_2$ are
eigenvectors of $u$. Hence, $u^*(\omega_1)$ must be $\pm \omega_1$
and $u^*(\omega_2)=\mp \omega_2$. Therefore
$u^*(x)=u^*(\omega_1/\omega_2)=  -x$ and $u^*(y)$ is $y$ or $-y$
depending on whether $u^*(\omega_2)$ is $-\omega_2$ or not.      In
any case, $y^2=P(x^2)$ for a  degree three  polynomial $P\in\Q[x]$,
and  the automorphism group of the curve is the  Klein group
generated by $u$ and the hyperelliptic involution $w$. {{ Moreover,
$X_0^*(N)$ has two bielliptic involutions $u$ and $u\cdot w$ and
both elliptic curves are bielliptic quotients.}}
\end{rem}
\vskip 0.2 cm Now, we will apply two sieves to discard some values
of $N$. Both are based on the values of $|X_0^*(N)(\F_{p^n})|$  for
a prime $p\nmid N$. The first of them uses \cite[Theorem 2.1]{Go17},
which  allows us to detect   some curves $X/\Q$ without involutions
defined over $\Q$,{{ because $\Aut_\Q(X)\hookrightarrow
\Aut_{\F_p}(X/\F_p)$ for a prime $p$ of good reduction for  $X$ (see
\cite[Prop.10.3.38]{LiuBook})}}. More precisely, for such  a prime
$p$  and an integer $n\geq 1$, consider the sequence
$$P_p(n):= \operatorname{mod\,} [(\sum_{d|n}\mu(n/d)|X(\F_{p^n})|)/n,2]\, $$
 where  $\operatorname{mod\,} [r,2]$ denotes $0$ or $1$ depending on whether $r$ is even or not, and
$\mu$ is the {{Moebius}} function. Set $Q_p(2 k+1)=\sum_{n\geq
0}^k(2 n+1) P_p(2n+1)$. If $X_0^*(N)$  has an involution defined
over $\Q$, then
 $$Q_p(2k+1)\leq 2 g^*_N+2  \text{ for all $k\geq 0$}\,.$$

\begin{prop}\label{G}\label{prop2}
The curve $X_0^*(N)$ is not bielliptic and, moreover,
$\Aut(X_0^*(N))$ is the trivial group for the following values of $N$:
$$\begin{array}{c}259,267,301,305,393,395,427,445,581,795,903,1001,1015\,.\end{array}$$
\end{prop}

\noindent{\bf Proof.}
For  $f\in\New_M$, denote by $K$ the number field $\Q(\{a_n(f)\}_{n\geq 1})$. Let $p$ be a prime not dividing $M$.  By the Eichler-Shimura congruence, the characteristic polynomial of $\Frob_p$ acting on the Tate module of $A_f$ is
$$
\prod_{\sigma \colon K\hookrightarrow \overline{\Q}}(x^2-a_p(f^\sigma)+p)\,,
$$
where $\sigma$ runs over the set  of all $\Q$-embeddings of $K$ into a fixed algebraic closure of $\Q$. The jacobian of $X_0^*(N)$ is isogenous over $\Q$ to the product $\prod_{1\leq M|N}\prod_{f\in\New_M^*/G_{\Q}}A_f$, where $G_{\Q}$ denotes the absolute Galois group  $\Gal(\overline{\Q}/\Q)$.

To compute $|X_0^*(N)(\F_{p^n})|$, we proceed as follows.
By using Magma,  we determine $\New_M^*$ for all $ M|N$ and, then, the characteristic polynomial $R_p(x)$ of $\Frob_p$ acting on the Tate module of $\Jac(X_0^*(N))$ is obtained as follows
$$
R_p(x)= \prod_{M|N} \prod_{f\in{\New_M}^*}(x^2-a_p(f) \,x+p)=\prod_{i=1}^{2\,g_N^*}(x-\alpha_i)\,.
$$
Finally,
$$
|X_0^*(N)(\F_{p^n})|=p^n+1-\sum_{i=1}^{2\,g_N^*} \alpha_i^n\,.
$$
The statement follows from these computations:
   $$
\begin{array}{|c|c|r|}
\hline N  & Q_p(2k+1)  & 2 g_N^* +2 \\ \hline\hline 259=7\cdot 37 &
{{  Q_2(9)=17}}        & 10\\ \hline 267=3\cdot 89 & {{Q_2(7)=15}} & 10\\
\hline
 301=7\cdot 43 &  {{Q_2(11)=22}}          & 14\\ \hline
305=5\cdot 61 &Q_2(7)=15  & 10 \\ \hline 393=3\cdot 131 &
{{Q_2(11)=23}} &12\\ \hline 395=5\cdot 79 & Q_2(7)=15& 10\\
\hline 427=7\cdot 61 &  Q_2(11)=27           & 18\\ \hline
445=5\cdot 89& {{Q_3(9)=22}} & 16\\  \hline
         581= 7\cdot 83  &Q_2(13)=20 &18\\ \hline
\end{array}
\quad
        \begin{array}{|l|c|c|}
         \hline
         N & Q_p(2k+1)& 2g_N^*+2\\ \hline
                \hline
   \hline 795=3\cdot5\cdot 53&Q_2({{13}})=27 &22 \\
          \hline 903=3\cdot 7\cdot 43&{{Q_2(13)= 31}}& 24\\
             \hline 1001=7\cdot 11\cdot 13&{{Q_2(11)=22}} & 18\\
        \hline 1015=5\cdot 7\cdot 29&Q_2(13)=30 & 24\\
                  \hline
                      \end{array}\quad \Box
$$

\vskip 0.2 cm
The second sieve  is based on the following fact. For a degree two morphism  of  curves  $X \rightarrow Y$ defined over $\Q$
 and a prime $p$  of good
reduction for $X$, one has
$$
       |X(\F_{p^k}) |-2 |Y(\F_{p^k})|     \leq 0\,, \text{ for all }k>0\,.
$$
\begin{prop}
   The pairs $(N,E)$ in the set
   $$\begin{array}{c}
   \{{{(273,91a)}},(385,{{77a}}),(415,83a),{{(429,143a)}}, {{(435,145a)}},
   {{(455,91a)}},\\(465,155c),
   (555,185c),(615,123b),(705,141d),
   (715,65a),(715,143a),\\{{(861,123b)}},(987,141d),
   (1155,77a),(1365,65a),(1365,455a)\}
   \end{array}$$   are            not bielliptic.
    In particular,  the curve $X_0^*(N)$ is         not bielliptic for the following values of $N$:
    $${{273}},385,415,{{429}},{{435}}, 465, {{555}}, 615,705, 715, {{861}},987,1155,1365\,.$$
             .
\end{prop}
\noindent{\bf Proof.}
      $$ \begin{array}{|c|c|c|c|}       \hline
           N & p^k & E& |X_0^* (\F_{p^k})|-     2 |E(\F_{p^k})| \\
\hline {{273}}&8&91a&3\\
           \hline
           385 &  9  & {{77 a}}  &   {{ 3}}        \\
           \hline 415&9&83a& 4\\
           \hline {{429}}&16&143a&5\\
           \hline {{455}}&4&91a&1\\
           \hline {{435}}&8&145a&2\\
           \hline 465&7&155c&2\\
           \hline 555&2&185c&1\\
           \hline 615&11&123b&4\\
           \hline 705&4&141d&1\\
           \hline 715&9&65a&7\\
           \hline 715&9&143a&1\\
           \hline {{861}}&4&123b&1\\
           \hline 987&25&141d&10\\
           \hline 1155&2&77a&1\\
           \hline 1365&4&65a&1\\
           \hline 1365&2&455a&3\\

           \hline
      \end{array} \quad \Box
      $$
      \vskip 0.2 cm
    After applying the two sieves,     the following possibilities for the pairs $(N,E)$,  ordered by the genus, remain.

    $$
    \begin{array}{|l|r|c|}
     \hline
     N & g_N^*& E\\ \hline
     183& 3& 61 a\\
     \hline185& 3&37a\\
     \hline 249&3&83a\\
      &&249b\\
      \hline 303&3&101 a\\
     \hline 455&3&65a\\\hline
\end{array}
\quad
\begin{array}{|l|r|c|}
     \hline
     N & g_N^*& E\\ \hline
      \hline 399&4&57a\\
\hline
     \hline 237&5&79a\\


      \hline 645&5&129a\\
      &&215a\\
     \hline
     \end{array}
        $$

     $$\text{Table 2} $$

Finally, in order to decide which  values $N$ in Table 2 correspond to bielliptic curves, we shall use equations.

We recall that, for a nonhyperelliptic curve $X$ defined over $\C$
with genus $g>2$, the image of the canonical map $ X\rightarrow
\PP^{g-1}$ is the common zero locus of a set of homogeneous
polynomials of degree $2$  and $3$, when $g>3$,  or of a homogenous
polynomial of degree $4$, if $g=3$.

More precisely, assume that  $X$ is defined over $\Q$ and  choose a basis $\omega_1,\cdots\omega_g$ of $\Omega^1_{X/\Q}$.
For any integer $i\geq 2$, let us denote by $\cL_i$ the $\Q$-vector space of homogeneous polynomials $Q\in\Q [x_1,\cdots,x_g]$
of  degree $i$ that satisfy $Q(\omega_1,\cdots,\omega_g)=0$. Of course,  $\dim \cL_i\leq \dim \cL_{i+1}$ because one has $x_j\cdot Q \in\cL_{i+1}$ for all $Q\in
\cL_i$  and for $1\leq j\leq g$.

If $g=3$, then $\dim \cL_2=\dim \cL_3=0$ and $\dim \cL_4=1$. Any
generator of $\cL_4$ provides an equation for $X$. For $g>3$, $\dim
\cL_2=(g-2)(g-3)/2>0$ and a basis of $\cL_2\bigoplus\cL_3'$ provides
a system of equations for $X$, {{ where $\cL_3'$ is any complement
of the vector subspace of $\cL_3$ consisting of all polynomials that
are multiples of a polynomial in $\cL_2$}}. When $X$ is neither
trigonal nor a smooth plane quintic ($g=6$), it suffices to take a
basis of  $\cL_2$.

For the curve  $X_0^*(N)$ there exists a set of normalized
eigenforms $g_1,\cdots,g_k\in S_2(\Gamma_0(N))^{B(N)}$ such that
$\Jac(X_0^*(N))\stackrel{\Q}\sim A_{g_1}\times \cdots\times
A_{g_k}$, where  the symbol $\stackrel{\Q}\sim$ means isogenous over
$\Q$.  These abelian varieties are simple and pairwise nonisogenous
over $\Q$. Hence, any involution $u$ of the curve leaves stable
$A_{g_i}$ and  acts on $\Omega^1_{A_{g_i}}$ as the product by $-1$
or the identity, {{because the endomorphism algebra
$\End_{\Q}A_{g_i}\otimes \Q$ is isomorphic to a (totally real)
number field}}.

We choose a basis $\{\omega_1,\cdots,\omega_{g_N^*}\}$ of $\Omega^1_{X_0^*(N)/\Q}$ obtained as the union
of bases of all $\Omega^1_{A_{g_i}/\Q}$.  An involution $u$ of $X_0^*(N)$ induces a linear map $u^*:\Omega^1_{X_0^*(N)/\Q}\rightarrow \Omega^1_{X_0^*(N)/\Q}$ sending $(\omega_1,\cdots,\omega_{g_N^*})$ to $(\varepsilon_1\omega_1,\cdots,\varepsilon_n\omega_{g_N^*})$ with $\varepsilon_i=\pm 1$  for all $i\leq {g_N^*}$ and satisfying
\begin{equation}\label{involution}
Q(\varepsilon_1x_1,\cdots,\varepsilon_{g_N^*}  x_{g_N^*})\in \cL_i \text{ for all } Q\in\cL_i \text{ and for all } i\,.
\end{equation}
Conversely, for a linear map $u^*$ as above satisfying condition (\ref{involution}), only one of the two maps $\pm u^*$ comes from an involution of the curve, because we are assuming that $X$ is nonhyperelliptic.

 \vskip 0.2 cm
 We particularize this fact to our case.
\begin{lema}\label{con}
 Assume $X_0^*(N)$ is nonhyperelliptic. Let $\omega_1,\cdots,\omega_{g_N^*}$ be a basis of $\Omega^1_{X_0^*(N)/\Q}$ as above, such that
  $\omega_1$ is the differential attached to an  elliptic curve $E$. Then,  the pair $(N,E)$ is bielliptic if, and only, if
\begin{equation} \label{eq}Q(-x_1,x_2,\cdots,x_{g_N^*-1}, x_{g_N^*})\in \cL_i \text{ for all } Q\in\cL_i \text{ and for all } i\,.\end{equation}
\end{lema}

\noindent{\bf Proof.}
If $u$ is an involution of $X_0^*(N)$ such that $E$ is $\Q$-isogenous to $X_0^*(N)/\langle u\rangle$, then  $u^*(\omega_1)=\omega_1$ and $u^*(\omega_i)=-\omega_i$ for $i>1$. Hence,  condition (\ref{eq})  is satisfied. Conversely, since the curve is nonhyperelliptic
the  condition (\ref{eq}) implies that only one of the two linear maps
$$
(\omega_1,\omega_2,\cdots ,\omega_{g_N^*-1}, \omega_{g_N^*})\mapsto\pm (-\omega_1,\omega_2,\cdots ,\omega_{g_N^*-1},\omega_{g_N^*})
$$
comes from  an involution $u$ of the curve.  The genus $g_u$  of the
curve  $X_0^*(N)/\langle u\rangle$ agrees with the number of
differentials  $\omega_i$ invariant under the action of  $u$. When
$g_N^*>3$,  it follows that $g_u$  must be $1$ because it cannot be
$g_N^*-1$ due to  Riemann-Hurwitz formula. For $g_N^*=3$,  the genus
$g_u$ must be different from $2$, since otherwise  the curve would
be  hyperelliptic (cf. \cite[Lemma 5.10]{Accola}). \hfill $\Box$

\begin{rem}
When $g_N^*>4$, $\dim \cL_2 >1$. If $\omega_j$ is the differential attached to an elliptic curve, we need to check that the vector space
$$
\cL_{2,j}:=\{Q\in\cL_2\colon Q(x_1,\cdots, -x_j,\cdots,
x_{g_N^*})\in\cL_2\}
$$
is $\cL_2$. Note that
$$\cL_{2,j}=\{Q\in\cL_2\colon Q(x_1, \cdots, x_j,\cdots, x_{g_N^*})=Q(x_1,\cdots, -x_j,\cdots, x_{g_N^*})\}
$$
Indeed, if $Q\in \cL_{2,j}$, then $H:=Q(x_1, \cdots, x_j,\cdots,
x_{g_N^*})-Q(x_1,\cdots, -x_j,\cdots x_{g_N^*})\in\cL_2$. Therefore,
$H=x_j P$ for an homogenous polynomial $P\in
\Q[x_1,\cdots,x_{g_N^*}]$ of degree at most $1$. Hence, $P$ must be
$0$,  otherwise $P(\omega_1,\cdots,\omega_{g_N^*})=0$.
\end {rem}

{{
\begin{rem}
Recall that,  for each one of the normalized  eigenforms  $g_i\in
S_2(\Gamma_0(N))^{B(N)}$,  there is $f_i\in\New_M^*$ such that
$g_i=\sum_{d|N/M} d \,f_i(q^d)$ and $A_{g_i}\stackrel{\Q}\sim
A_{f_i}$. To get a basis of $\Omega^1_{A_{g_i}/\Q}$ we can proceed
as follows. If $\dim A_{f_i}=1$, we take as basis $g_i(q)\, dq/q$.
When $\dim A_{f_i}=r>1$,  the endomorphism algebra $\End_\Q
(A_{f_i})\otimes \Q$ is generated by Hecke  operators and  is
isomorphic to the totally real number field
$K_i:=\Q(\{a_n(f_i)\}_{n>0})$ of degree $r$. Let $\cI$ be the set of
$\Q$-embeddings of $K_i$ into a fixed algebraic closure of $\Q$.
  For every $a\in K_i$ there is a Hecke operator $T$ such that $T(f_i^{\sigma})=a^\sigma f_i^\sigma$ for all $\sigma\in\cI$. The two cusp-form $h=\sum_{\sigma\in\cI} f_i^\sigma$ is nonzero because the coefficient of $q$ is $g$. Hence,  taking $T$ such that $a$ is a primitive element of $K_i$,  the set  $$f'_j=T^j (h)=\sum_{\sigma\in\cI} (a^j )^\sigma  f_i ^\sigma\in\Q [[q]],\,\, 0\leq j\leq r-1\,,$$
is a basis of the vector space  spanned by $f_i^\sigma$.
Therefore,
$$\omega_j'=(\sum_{d|N/M} d \,f_i'(q^d))\, dq/q\,\,\,  0\leq j\leq r-1\, $$
is a basis of $\Omega^1_{A_{g_i}/\Q}$. One can  take $a$ as
the value provided by Magma in the $q$-expansion of $f_i$ and, in this
case, $f_j'\in\Z [[q]]$. The curve  $X_0^*(N)$ is determined by the first Fourier coefficients of the chosen basis for $\Omega^1_{X_0^*(N)/\Q}$ (cf \cite[Proposition 2.8]{BGGP}).  In order to get shorter equations, it is suitable to replace the basis $f_j'$ with a basis
of the $\Z$-module $\left(\bigoplus_{i=1}^g\Q f_i'\right)\cap \Z[[ q]]$.

\end{rem} }}

\vskip 0.2 cm

\begin{prop}\label{prop4}
Among the curves of genus three $X_0^*(183 )$, $X_0^*(185 )$,
$X_0^*( 249)$, $X_0^*(303)$, and $X_0^*(455 )$,  only $X_0^*(183 )$,
$X_0^*( 249)$, $X_0^*(303)$  and $X_0^*(455 )$  are bielliptic. The
corresponding elliptic quotients are labeled as $E61a1$, $E249b1$,
$E101a1$ and $E65a1$, respectively. In all these cases, the
automorphism group has order 2. The automorphism groups of the
remaining curves are trivial.
\end{prop}

\noindent{\bf Proof.}  For these values of $N$, the splitting of the jacobian of $X_0^*(N)$, $J_0^*(N)$, is as follows:
$$\begin{array}{rrcrr}
J_0^*(183)\stackrel{\Q}\sim &\prod_{i=1}^2A_{f_i}\,, &A_{f_1}\stackrel{\Q}\sim E61a\,,\,& f_2\in\New_{183}^*\,,\,&\dim A_{f_2}=2\,,\\
J_0^*(185)\stackrel{\Q}\sim &\prod_{i=1}^3A_{f_i}\,, &A_{f_1}\stackrel{\Q}\sim E37a\,,\, &A_{f_2}\stackrel{\Q}\sim E185{{a}}\,, &A_{f_3}\stackrel{\Q}\sim E185c\,,\\
J_0^*(249)\stackrel{\Q}\sim &\prod_{i=1}^3A_{f_i}\,,\, &A_{f_1}\stackrel{\Q}\sim E83a\,,&\,A_{f_2}\stackrel{\Q}\sim E249a\,,\, &A_{f_3}\stackrel{\Q}\sim E249b\,,\\
J_0^*(303)\stackrel{\Q}\sim &\prod_{i=1}^2A_{f_i}\,, &A_{f_1}\stackrel{\Q}\sim E101a\,,\, &f_2\in\New_{303}^*\,,\, &\dim A_{f_2}=2\,,\\
J_0^*(455)\stackrel{\Q}\sim &\prod_{i=1}^3A_{f_i}\,, &A_{f_1}\stackrel{\Q}\sim E65a\,,\,&A_{f_2}\stackrel{\Q}\sim {{E91a}}\,,\,&  A_{f_3}\stackrel{\Q}
\sim {{E455a}}\,.\\
\end{array}
$$
We take a basis of $\Omega^1_{X_0^*(N)/\Q}$ following the order exhibited in the splitting of its jacobian and we obtain the following generators  $Q\in\cL_4$:
$$\begin{array}{c|r}
N & Q \phantom{ccccccccccccccccccccccccccccccccccccccc}\\ \hline
183 &x^4 - 10 x^2 y^2 + 9 y^4 - 24 x^2 y z + 24 y^3 z + 32 y^2 z^2 +
 32 y z^3 - 16 z^4\\
 185&
2816 x^4 + 768 x^3 y + 1728 x^2 y^2 - 243 y^4 - 5888 x^3 z +
 1152 x^2 y z + 1728 x y^2 z \\& + 324 y^3 z + 192 x^2 z^2 -
 1152 x y z^2 + 918 y^2 z^2 - 2624 x z^3 + 852 y z^3 - 571 z^4\\
  249 &16 x^4 - 64 x^3 y - 12 x^2 y^2 + 44 x y^3 + 97 y^4 - 180 x^2 z^2 +
 36 x y z^2 - 18 y^2 z^2 + 81 z^4\\
 303 &x^4 + 2 x^2 y^2 - 3 y^4 - 16 y^2 z^2 + 16 y z^3 - 16 z^4\\
 455 &81 x^4 - 162 x^2 y^2 - 79 y^4 - 324 x^2 y z + 244 y^3 z +
 192 y^2 z^2 + 64 y z^3 - 16 z^4\\
 \end{array}
 $$
By Lemma \ref{con}, only the curves corresponding to $N=183,
249,303, 455$ are bielliptic with an only bielliptic involution $u$. The affine equations for the bielliptic quotients are
$$\begin{array}{c|r}
N &  X_0^*(N)/\langle u\rangle\phantom{cccccccccccccccccccccccccccccc}\\ \hline
183 &x^2 - 10 x  + 9  - 24 x z + 24  z + 32  z^2 +
 32  z^3 - 16 z^4=0\\
 249 &97 + 44 x - 12 x^2 - 64 x^3 + 16 x^4 - 18 z + 36 x z -
 180 x^2 z + 81 z^2=0\\
 303 &x^2 + 2 x  - 3  - 16 z^2 + 16 z^3 - 16 z^4=0\\
  455 &81 x^2 - 162 x  - 79  - 324 x  z + 244  z +
 192  z^2 + 64  z^3 - 16 z^4=0\\
 \end{array}
 $$
which  have genus one and their $j$-invariants are $-\frac{912673}{61}$,  $\frac{357911}{249}$, $\frac{262144}{101}$ and $\frac{117649}{65}$.
They correspond to the elliptic curves $E61a1$, $E249b1$, $E101a1$ and $E65a1$.

Taking into account the splitting of their jacobians and their
equations, all their automorphism groups have order $2$. The
remaining curves have trivial automorphism groups. For instance, for
$N=249$, the linear map $(\omega_1,\omega_2,\omega_3)\mapsto
(\omega_1,\omega_2,-\omega_3)$ is the only option to be considered
and $Q(-x,y,z)\notin\cL_4$. For $N=185$, none of the polynomials
$Q(-x,y,z)$,$Q(x,-y,z)$, $Q(x,y,-z)$ lies in  $\cL_4$. \hfill $\Box$

\begin{prop}\label{prop5} The curve of genus four  $X_0^*(399)$ is not bielliptic and its automorphism group is trivial.

\end{prop}
\noindent{\bf Proof.} The splitting of  $J_0^*(399)$ is:
$$\begin{array}{rrrrrr}
J_0^*(399)\stackrel{\Q}\sim &A_{f_1}\times A_{f_2}\times A_{f_3}\,, &A_{f_1}\stackrel{\Q}\sim E57a\,,&\,A_{f_2}\stackrel{\Q}\sim E399a\,,\,\, &f_3\in\New_{133}^*\,,\, &\dim A_{f_3}=2\,.
\end{array}
$$
In this case $\dim \cL_2=1$. As in the previous  proposition, we take a basis of $\Omega^1_{X_0^*(399)\Q}$ following the order exhibited in the splitting of its jacobian. Next, we show a generator $Q(x,y,z,t)\in \cL_2$:
$$\begin{array}{c|r}
N & Q \phantom{ccccccccccccccccccccccccccccccccccccccc}\\ \hline
 399 & -99 t^2 + 90 t x + 125 x^2 + 189 t y + 80 x y - 151 y^2 + 306 t z -
 105 x z + 42 y z + 9 z^2 \\
 \end{array}
 $$
 Since $Q(-x,y,z,t)\notin \cL_2$, the curve is not bielliptic. The conditions $Q(x,-y,z,t)$,$Q(-x,-y,z,t)$ $\notin\cL_2$ imply that the curve does  not have any nontrivial involutions.
\hfill $\Box$

\begin{prop} The curves of genus five $X_0^*(237)$ and  $X_0^*(645)$ are not bielliptic.

\end{prop}
\noindent{\bf Proof.}
The splitting of  $J_0^*(N)$ is:
$$\begin{array}{rrrrrr}
J_0^*(237)\stackrel{\Q}\sim& \prod_{i=1}^2A_{f_i}\,, & &A_{f_1}\stackrel{\Q}\sim E79a\,,\, &f_2\in\New_{237}^*\,,\,&\dim A_{f_2}=4\,,\\
J_0^*(645)\stackrel{\Q}\sim &\prod_{i=1}^4A_{f_i}\,, &A_{f_1}\stackrel{\Q}\sim E43a\,,\,&A_{f_2}\stackrel{\Q}
\sim E129a\,,&\,A_{f_3}\stackrel{\Q}\sim E215a\,,\,&f_4\in\New_{645}^*\,,\,\\& && &
  &\dim A_{f_4}=2\,.
\end{array}
$$
Now, $\dim \cL_2=3$. For $N=237$,  $\dim \cL_{2,1}=0$ and  for $N=645$ we also have  $\dim\cL_{2,2} =\dim \cL_{2,3}=0$.
 \hfill $\Box$

\vskip 0.2 cm
As a consequence, we obtain the statement of Theorem \ref{main} for $N$ odd.
\begin{cor}
For $N$ odd, $X_0^*(N)$ is bielliptic if, and only if, $N\in
\{129,183,215,249,303, 455 \}$. For these values of $N$, the automorphism group has order
$2$ when $g_N^*>2$, otherwise it is the Klein group.
\end{cor}

\section{Even case}
  By applying Lemma \ref{psi}, we determine a finite set of possible values of $N$.
  Then, we proceed as in the odd case and
we obtain the pairs $(N,E)$ exhibited in Table 3 together the genera
of $X_0(N)$ and $X_0^*(N)$. { {As in the odd case, for $g_N^*=2$ we
can discard $N=154,286$ because both curves $X_0^*(N)$ have an
elliptic quotient of conductor $N$ that does not satisfy Lemma
\ref{degree}}}. Table 3  is divided into  4
 cases:  $3\nmid N$,  $6\mid N$ and $30\nmid N$,  $30\mid N$ and $210\nmid N$ and, finally,   $210\mid N$.

\begin{center}
\begin{tiny}
$$
\begin{array}{|l|r|r|r|c|}
\hline
N , 3\nmid N&g_N&g_N^*&M &\text{Label $E$}\\
        \hline
  106=2\cdot 53& 12& 2&106 &b\\
          &   &   &   53 & a\\ \hline
  122=2\cdot 61& 14& 2& 122 &a\\
          & &&61&a\\
  \hline 158=2\cdot 79& 19& 2&158 &b\\
              &  &   & 79&a\\
  \hline 166=2\cdot 83& 20& 2&166&a\\
    &  &  &83 &a\\
  \hline 178=2\cdot 89&21&3&89&a\\
  \hline 202=2\cdot 101&24&4&101&a\\
  \hline 262=2\cdot 131& 32& 4&131&a\\
\hline 394=2\cdot 197& 48& 10&197&a\\
\hline 458=2\cdot 229& 56& 10&229&a\\
\hline
\hline 290= 2\cdot 5\cdot 29 & 41& 3&58&a\\
 &&&145&a\\
\hline 310= 2\cdot 5\cdot 31 & 45& 3&155&c\\
\hline 370=  2\cdot5\cdot 37 & 53& 4&
185&a\\
&&&185&c\\
&&&370&a\\
\hline 410= 2\cdot 5\cdot 41 & 59& 5&82&a\\
\hline 430= 2\cdot 5\cdot 43 & 63& 3&43&a\\
&&& 215&a\\
\hline 530= 2\cdot 5\cdot 53 & 77& 7&
106&b\\
&&&265&a\\
\hline 574= 2\cdot 7\cdot 41 & 81& 5&82&a\\
\hline 590=  2\cdot5\cdot 59 &87& 6&118&a\\
\hline 602= 2\cdot 7\cdot 43 & 85& 9&43&a\\
\hline 710=2\cdot 5\cdot 71 & 105& 7&142&b\\
\hline 742= 2\cdot 7\cdot 53 & 105& 10&
371&a\\
\hline 770= 2\cdot 5\cdot 7\cdot 11 & 137& 5&77&a\\
&&&154&a\\
\hline 910= 2\cdot 5\cdot 7\cdot 13 & 161& 5&65&a\\
&&&91&a\\
&&&455&a\\
\hline 1190=2\cdot 5\cdot 7\cdot 17 & 209& 8&238&b\\
\hline
\end{array}  \quad
\begin{array}{|l|r|r|r|c|}
\hline
6|N,30\nmid N &g_N&g_N^*&M &\text{Label $E$}\\
\hline
\hline 246= 2\cdot 3\cdot 41 & 39& 3&82&a\\
&&&123&b\\
\hline 258=  2\cdot 3\cdot 43 & 41& 3&43&a\\
&&&129&a\\
\hline 282= 2\cdot 3\cdot 47 & 45& 3&141&d\\
\hline 318= 2\cdot 3\cdot 53 & 51& 3&53&a\\
&&&106&b\\
\hline 354= 2\cdot 3\cdot 59 & 57& 4&118&a\\
\hline 366=  2\cdot 3\cdot 61 & 59& 4&61&a\\
&&&122&a\\
\hline 402= 2\cdot 3\cdot 67 & 65& 5& 201&a\\
&&&201&c\\
\hline 426= 2\cdot 3\cdot 71 & 69& 4&142&b\\
\hline 438= 2\cdot 3\cdot 73 & 71& 5&219&a\\
&&&219&c\\
\hline 474= 2\cdot 3\cdot 79 & 77& 7&79&a\\
&&&158&b\\
\hline 498= 2\cdot 3\cdot 83 & 81& 7&83&a\\
&&&166&a\\
&&&249&a\\
&&&249&b\\
\hline 534= 2\cdot 3\cdot 89 & 87& 8&89&a\\
\hline 582= 2\cdot 3\cdot 97 & 95& 7&291&c\\
\hline 606= 2\cdot 3\cdot 101 & 99& 8&101&a\\
\hline 642= 2\cdot 3\cdot 107 & 105& 9&214&b\\
\hline 786=  2\cdot 3\cdot 131& 129& 11&131&a\\
\hline
\hline 462= 2\cdot 3\cdot 7\cdot 11 & 89& 3&77&a\\
&&&154&a\\
\hline 546= 2\cdot 3\cdot 7\cdot 13 & 105& 4&91&a\\
\hline 714= 2\cdot 3\cdot 7\cdot 17 & 137& 5&102&a\\
&&&238&b\\
\hline 798= 2\cdot 3\cdot 7\cdot 19 & 153& 5&57&a\\
&&&399&a\\
\hline 858=  2\cdot 3\cdot 11\cdot 13 & 161& 6&143&a\\
&&&286&c\\
\hline 966= 2\cdot 3\cdot 7\cdot 23 & 185& 8&138&a\\
\hline 1122= 2\cdot 3\cdot 11\cdot 17 & 209& 9&
374&a\\
\hline 1254= 2\cdot 3\cdot 11\cdot 19 & 233& 12&57&a\\
\hline
\end{array}
$$
$$
\begin{array}{|l|r|r|r|c|} \hline
30|N, 210\nmid N &g_N&g_N^*&M &\text{Label $E$}\\
\hline
\hline 390=  2\cdot 3\cdot 5\cdot 13 & 77& 2&65&a\\
&&&390&a\\
\hline 510= 2\cdot 3\cdot 5\cdot 17 & 101& 3&102&a\\
\hline 570= 2\cdot 3\cdot 5\cdot 19 & 113& 4&57&a\\&&&190&b\\&&&285&b\\
\hline 690=  2\cdot 3\cdot 5\cdot 23 & 137& 6&138&a\\
\hline 870= 2\cdot 3\cdot 5\cdot 29 & 173& 7&58&a\\&&&145&a\\&&&290&a\\
\hline 930= 2\cdot 3\cdot 5\cdot 31 & 185& 8&155&c\\
\hline 1110=  2\cdot 3\cdot 5\cdot 37 & 221& 9&
185&c\\
\hline
\end{array} \quad \begin{array}{|l|r|r|r|c|}
\hline 30|N, 210\nmid N&g_N&g_N^*&M&\text{Label $E$}\\
\hline 1230= 2\cdot 3\cdot 5\cdot 41 & 245& 10&
123&b\\
&&&{{615}}&a\\
\hline 1290= 2\cdot 3\cdot 5\cdot 43 & 257& 10&
129&a\\
&&&215&a\\
\hline 1410= 2\cdot 3\cdot 5\cdot 47 & 281& 14&141&d\\
&&&705&a\\
\hline 1590= 2\cdot 3\cdot 5\cdot 53 & 317& 14&53&a\\&&&265&a\\&&&795&a\\
\hline \hline
210|N &g_N&g_N^*&M &\text{Label $E$}\\
\hline

\hline 2310= 2\cdot 3\cdot 5\cdot 7\cdot 11 & 561& 12&77&a\\&&&1155&a\\
\hline 2730= 2\cdot 3\cdot 5\cdot 7\cdot 13 & 657& 14&65&a\\
&&&455&a\\
\hline
\end{array}
$$

\end{tiny}
\end{center}
$$
    \text{Table 3}
$$

As in the odd case, first we examine the hyperelliptic curves.
\begin{prop} {{ All the curves of genus two appearing in Table 3, i.e. $X_0^*(106)$, $X_0^*(122)$,  $X_0^*(158)$, $X_0^*(166)$ and $X_0^*(390)$,
 are bielliptic}}. \end{prop}
\noindent{\bf Proof.} For $N\in\{106,122, 154, 158,166, 286, 390\}$,
the jacobian of $X_0^*(N)$ is isogenous over $\Q$ to the product of
two elliptic curves $E$ and $F$ of conductors $N$ and $M<N$,
respectively.
We have that $M=N/6$ for $N=390$ and $M=N/2$
otherwise. Let $f_E\in\New_N^*$ and $f_F\in\New_M^*$ be the
corresponding newforms attached to these elliptic curves. The
functions $$ x:= \frac{f_E(q)}{\sum_{d|N/M} d\,f_F(q^d)}\,, \quad
y:= \frac{ 2\,q\,{\rm d\, }x }{(\sum_{d|N/M} d\,f_F(q^d)){\rm d\, } q}\,,
$$
  provide the following equations
$$
\begin{array}{rccr}
X_0^*(106)\colon &  y^2 &=&  4 x^6+17 x^4-6x^2+1\,,  \\
X_0^*(122)\colon&  y^2 &=& 4 x^6  + x^4+10 x^2+1\,,\\
X_0^*(158)\colon& y^2 &=&  - 2 x^6+11 x^4+8 x^2-1\,, \\
X_0^*(166)\colon &y^2&=&- 4 x^6+17 x^4+2 x^2+1\,, \\
X_0^*(390)\colon & y^2&=&-(3 x^2+1) ( 4 x^4-7 x^2-1)\,.
\end{array}
$$
In all cases, one has the involutions $(x,y)\mapsto  \pm (-x,\pm y)$. \hfill $\Box$

\vskip 0.2 cm
Next, we  use the sieve based on
 \cite[Theorem 2.1]{Go17}.

\begin{prop}\label{prop8}
The curve $X_0^*(N)$ is not bielliptic and its automorphism group is trivial for the following values of $N$:
$$\begin{array}{c}394,458,582,602,710,786\,.\end{array}$$
\end{prop}

\noindent{\bf Proof.} For a prime $p\nmid N$, let $Q_p(2 k+1)$ be as in Proposition \ref{G}.
After the following computations,
   $$
\begin{array}{|c|c|r|} \hline
N  & Q_p(2k+1)  & 2 g_N^* +2 \\
\hline 394=2\cdot 197 & Q_3(9)=25 & 22\\ \hline
 458=2\cdot 229 &  Q_5(15)=36          & 22\\
\hline 582=2\cdot 3\cdot 97 & Q_{13}(9)=18 &16\\
\hline 602=2\cdot 7\cdot 43 & Q_5(11)=31& 20\\
        \hline 710=2\cdot5\cdot 71&Q_3(9)=25 &16 \\
             \hline 786=2\cdot 3\cdot 131&Q_5(11)=27&24\\
%
          \hline
            \end{array}
$$
the statement follows. \hfill $\Box$

\vskip 0.2 cm
Now, we apply the  sieve  based on the values of
$|X_0^*(N)(\F_{p^n})|-2|E(\F_{p^n})|$ and a modification for primes
$p$ dividing the conductor $N$.

\begin{prop}
   The pairs $(N,E)$ in the set
   $$
   \begin{array}{cccccc}
   \{(290,58a),& (370,185a),& (402,201a),& (410,82a),& (438,219c),&
    (474, 79a),\\   (474,158b),& (498, 83a),& (498, 166a),& (498, 249a),& (498, 249b),
    &(530,106b),\\ (530,265a),&    (534,89a), &   (574,82a),& (590, 118a),& (606,101a),
    & (642,214b)\\ (714,102a),&    (742,371a), &(770,77a),&(770,154a),&(798,57a),
    & (870,58a)\\ (910,65a),&     (930,155c),&
      (966, 138a),&   (1110,185c),&    (1122,374a),
           &(1190,238b),\\
           (1230,123b),& (1230,615a),& (1254,57a),&
           (1290,129a),&(1290,215a),&(1410,141d),\\
 (1410, 705a), &   (1590,53a),&             (1590,265a),&
           (1590, 795a),&     (2310,77a),& (2310,1155a),\\
                  (2730,65a), &             (2730,455a)\}
   \end{array}
   $$
are  not bielliptic.
    In particular, the curve $X_0^*(N)$ is not bielliptic for the following values of $N$:
    $$\begin{array}{l}410,474, 498, 530,534, 574, 590, 606, 642, 742, 770,930, 966,1110, 1122, 1190,
 \\ 1230, 1254, 1290, 1410,1590, 2310, 2730\,.
 \end{array}$$

\end{prop}
\noindent{\bf Proof.} We put $n(N,E,p^k)=|X_0^*(N) (\F_{p^k})|-
2 |E(\F_{p^k})|$.
\begin{tiny}
      $$ \begin{array}{|c|c|c|c|}       \hline
           N & p^k & E& n(N,E,p^k) \\
           \hline 290&9&58a&5\\
           \hline 370
              &3&185a&2\\
           \hline 402&5&201a&2\\
           \hline 410&13&82a&1\\
           \hline 438&25&219c&12\\
           \hline 474&25& 79a&10\\
                   &11&158b&7\\
           \hline 498&11&83a&5\\
                   &7& 166a&3 \\
             &7& 249a&1\\
            &13& 249b&5\\
           \hline 530&
            25&106b&35\\
            &7&265a&8\\
           \hline 534&5&89a&2\\
           \hline 574&9&82a&4\\
           \hline 590&9&118a&2\\

            \hline
           \end{array} \quad
           \begin{array}{|c|c|c|c|}       \hline
           N & p^k & E& n(N,E,p^k) \\
                 \hline 606&5&101a&2\\
            \hline 642&7&214b&9\\
           \hline 714&25&102a&6\\
                     \hline 742&
              9&371a&4\\
           \hline 770&9&77a&9\\
        &3&154a&1\\
            \hline 798&25&57a&3\\
           \hline 870&13&58a&1\\
           \hline 910&9&65a&3\\
           \hline 930&7&155c&2\\
           \hline 966&25&138a&7\\

                     \hline 1122&
           5&374a&5\\
           \hline 1190&3&238b&1\\
           \hline
           \end{array} \quad
           \begin{array}{|c|c|c|c|}\hline
            N & p^k & E&  n(N,E,p^k)\\
           \hline 1230&
        49&123b&19\\
           &11&615a&2\\
            \hline 1254&5&57a&1\\
           \hline 1290&
           49&129a&3\\
           &49&215a&3\\

           \hline 1410&49&141d&20\\
           &7&705a&10\\
           \hline 1590&49&53a&54\\
           &7&265a&6\\
           &7&795a&10\\
           \hline 2310&17&77a&1\\
            &17&1155a&9\\
           \hline 2730&11&65a&7\\
          &11&455a&3\\
           \hline
      \end{array}
      $$
\end{tiny}

For the pair $(N,E)=(1110,185c)$, we proceed as follows.
     Suppose  that the pair $(N,E)$ is bielliptic.
      For  $p=2$, we know that $X_0(1110)$ modulo $2$ is the copy of two curves $X_0(555)/\F_2$, and the normalization of $X_0^*(1110)/\F_2$ is the
      curve  $X_0^*(555) /\F_2     $ (cf. \cite{FM}). Since $2$ does not divide
       the conductor of $E$, then  $E/\F_2$ is an elliptic curve that   is  the quotient curve of $X_0^*(555)/\F_2$ by an
       involution defined over $\F_2$. Therefore, $n(555, E, 2^k)\leq 0$. We get $n(555,E,2)=1$ and, thus, the pair $(1110,E)$ can be
       discarded.
    \hfill $\Box$

      \vskip 0.2 cm
  After applying the two sieves, the following possibilities for the pairs $(N,E)$,  ordered by the genus, remain:
    $$
    \begin{array}{|l|r|c|}
     \hline
     N & g_N^*& E\\ \hline
\hline 178&3&89a\\
\hline 246&3&82a,123b\\
\hline 258&3&43a,129a\\
\hline 282&3&141d\\
\hline 290&3&145a\\
\hline 310&3&155c\\
\hline 318&3&53a,106b\\
\hline 430&3&43a,215a\\
\hline 462&3&77a,154a\\
\hline 510&3&102a\\
 \hline
  \end{array} \quad
\begin{array}{|l|r|c|}
\hline N&g_N^*&E\\
\hline
\hline 202&4&101a\\
\hline 262&4&131a\\
\hline 354&4&118a\\
\hline 366&4&61a,122a\\
\hline 370&4&185c,370a\\
\hline 426&4&142b\\
\hline 546&4&91a\\
\hline 570&4&57a,190b,285b\\
\hline
\end{array}
 \quad
\begin{array}{|l|r|c|}
\hline N&g_N^*&E\\
\hline
\hline 402 & 5 &201c\\
\hline 438&5&219a\\
\hline 714&5&238b\\
\hline 798&5&399a\\
\hline 910&5&91a,455a\\
\hline\hline
690&6&138a\\
\hline 858&6&143a,286c\\
\hline
\hline 870&7&145a,290a\\
\hline
\end{array}
$$

$$\text{Table 4} $$

\begin{prop} \label{prop10}
Among the  curves of genus three $X_0^*(178)$, $X_0^*(246 )$,
$X_0^*( 258)$, $X_0^*(282)$, $X_0^*( 290)$, $X_0^*(310 )$,
$X_0^*(318)$, $X_0^*( 430)$, $X_0^*(462 )$ and $X_0^*(510 )$, only
$X_0^*(178 )$, $X_0^*(246 )$,  $X_0^*(258)$, $X_0^*(290 )$,
$X_0^*(318 )$,$X_0^*(430 )$    and $X_0^*(510 )$  are bielliptic.
The corresponding elliptic quotients are labeled as   $E89a1$,
$E82a1$, $E43a1$, $E145a1$, $E53a1$, $E43a1$ and $E102a1$,
respectively. In all these cases, the automorphism group {{ of
$X_0^*(N)$ }} has order 2.
The automorphism
groups of the remaining curves are trivial.
\end{prop}
\noindent{\bf Proof.}  For these values of $N$, the splitting of the jacobian of $X_0^*(N)$, $J_0^*(N)$, is as follows:
$$\begin{array}{rrccr}
J_0^*(178)\stackrel{\Q}\sim &\prod_{i=1}^2A_{f_i}\,, &A_{f_1}\stackrel{\Q}\sim E89a\,,\, &f_2\in\New_{178}^*\,,\,&\dim A_{f_2}=2\,,\\
J_0^*(246)\stackrel{\Q}\sim &\prod_{i=1}^3A_{f_i}\,, &A_{f_1}\stackrel{\Q}\sim E82a\,,&\,A_{f_2}\stackrel{\Q}\sim E123b\,,\,\,&A_{f_3}\stackrel{\Q}\sim E246d\,,\\
J_0^*(258)\stackrel{\Q}\sim &\prod_{i=1}^3A_{f_i}\,, &A_{f_1}\stackrel{\Q}\sim E43a\,,\,&A_{f_2}\stackrel{\Q}\sim E129a\,,\,\,&A_{f_3}\stackrel{\Q}\sim E258a\,,\\
J_0^*(282)\stackrel{\Q}\sim &\prod_{i=1}^2A_{f_i}\,, &A_{f_1}\stackrel{\Q}\sim E141d\,,\, &f_2\in\New_{282}^*\,,\, &\dim A_{f_2}=2\,,\\
J_0^*(290)\stackrel{\Q}\sim &\prod_{i=1}^3A_{f_i}\,, &A_{f_1}\stackrel{\Q}\sim E58a\,, &\, A_{f_2}\stackrel{\Q}\sim E145a\,,\,\,&A_{f_3}\stackrel{\Q}\sim E290a\,,\\
J_0^*(310)\stackrel{\Q}\sim& \prod_{i=1}^2A_{f_i}\,, &A_{f_1}\stackrel{\Q}\sim E155c\,,\,\, &f_2\in\New_{310}^*\,,\,&\dim A_{f_2}=2\,.\\
J_0^*(318)\stackrel{\Q}\sim &\prod_{i=1}^3A_{f_i}\,, &A_{f_1}\stackrel{\Q}\sim E53a\,,\, &A_{f_2}\stackrel{\Q}\sim E106b\,,\,\,&A_{f_3}\stackrel{\Q}\sim E318c\,,\\
J_0^*(430)\stackrel{\Q}\sim &\prod_{i=1}^3A_{f_i}\,, &A_{f_1}\stackrel{\Q}\sim E43a\,,\, & A_{f_2}\stackrel{\Q}\sim E215a\,,\,, &A_{f_3}\stackrel{\Q}\sim E430a\,,\\
J_0^*(462)\stackrel{\Q}\sim & \prod_{i=1}^3A_{f_i}\,, &A_{f_1}\stackrel{\Q}\sim E77a\,,\,& A_{f_2}\stackrel{\Q}\sim E154a\,,\,\,&A_{f_3}\stackrel{\Q}\sim E462a\,,\\
J_0^*(510)\stackrel{\Q}\sim &\prod_{i=1}^2A_{f_i}\,, &A_{f_1}\stackrel{\Q}\sim E102a\,,\,& f_2\in\New_{85}^*\,,\,& \dim A_{f_2}=2\,,\\
\end{array}
$$
We take a basis of $\Omega^1_{X_0^*(N)\Q}$ following the order exhibited in the splitting of the jacobian, and we obtain the following generators  $Q\in\cL_4$:
$$\begin{array}{c|r}
N & Q \phantom{ccccccccccccccccccccccccccccccccccccccc}\\ \hline
178 &x^4 - 2 x^2 y^2 + y^4 - 12 x^2 y z - 4 y^3 z - 4 x^2 z^2 +
 20 y^2 z^2 + 32 y z^3\\
 246 &81 x^4 - 16 y^4 + 324 x^2 y z - 32 y^3 z - 162 x^2 z^2 + 48 y^2 z^2 -
 260 y z^3 + 17 z^4  \\
 258& 81 x^4 - 90 x^2 y^2 + 41 y^4 + 288 x^2 y z - 176 y^3 z - 36 x^2 z^2 +
 84 y^2 z^2 - 128 y z^3 - 64 z^4\\
 282&x^4 - 4 x^3 y + 6 x^2 y^2 - 4 x y^3 + y^4 - 6 x^3 z - 18 x^2 y z -
 27 x y^2 z - 30 y^3 z\\ & - 9 x^2 z^2 + 45 x y z^2 - 117 y^2 z^2 +
 54 x z^3 + 108 y z^3\\
 290 &36 x^2 y^2 + 27 y^4 + 32 x^3 z + 36 x y^2 z + 12 x^2 z^2 -
 126 y^2 z^2 - 84 x z^3 + 67 z^4\\
 310& x^4 + 2 x^3 y - 3 x^2 y^2 - 4 x y^3 + 4 y^4 - 81 x y^2 z -
 81 x y z^2 - 81 y^2 z^2 - 54 x z^3 + 54 y z^3 + 81 z^4\\
  318&9 x^4 - 10 x^2 y^2 + y^4 - 28 x^2 y z + 12 y^3 z + 20 x^2 z^2 -
 4 y^2 z^2 - 32 y z^3 + 32 z^4\\
 430 & 81 x^4 + 54 x^2 y^2 - 7 y^4 - 432 x^2 y z - 128 y^3 z - 108 x^2 z^2 +
 444 y^2 z^2 + 64 y z^3 + 32 z^4 \\
 462& 128 x^4 - 320 x^3 y + 264 x^2 y^2 - 44 x y^3 - y^4 - 448 x^3 z +
 48 x^2 y z + 492 x y^2 z - 92 y^3 z \\&+ 840 x^2 z^2 + 12 x y z^2 +
 498 y^2 z^2 - 972 x z^3 - 972 y z^3 + 567 z^4\\
 510 & 3 x^4 - 4 x^2 y^2 + y^4 + 18 x^2 y z - 10 y^3 z + 14 x^2 z^2 +
 18 y^2 z^2 + 40 y z^3 + 16 z^4
 \end{array}
 $$
By Lemma \ref{con}, only the curves corresponding to $N=178, 246,
258, 290, 318, 430, 510$ are bielliptic  and  only have a bielliptic involution $u$. The affine equations for the bielliptic quotients are
$$\begin{array}{c|r}
N &  X_0^*(N)/\langle u\rangle\phantom{cccccccccccccccccccccccccccccc}\\ \hline
178 &-4 x + x^2 + 32 y - 12 x y + 20 y^2 - 2 x y^2 - 4 y^3 + y^4=0\\
 246 &17 - 162 x + 81 x^2 - 260 y + 324 x y + 48 y^2 - 32 y^3 - 16 y^4=0\\
  258&-64 - 36 x + 81 x^2 - 128 y + 288 x y + 84 y^2 - 90 x y^2 - 176 y^3 +
 41 y^4=0\\
  290 &36 x^2 y + 27 y^2 + 32 x^3  + 36 x y+ 12 x^2 -
 126 y  - 84 x + 67 =0\\
 318 & 32 + 20 x + 9 x^2 - 32 y - 28 x y - 4 y^2 - 10 x y^2 + 12 y^3 + y^4=0\\
  430& 32 - 108 x + 81 x^2 + 64 y - 432 x y + 444 y^2 + 54 x y^2 - 128 y^3 -
 7 y^4=0\\
 510 &16 + 14 x + 3 x^2 + 40 y + 18 x y + 18 y^2 - 4 x y^2 - 10 y^3 + y^4=0\\
 \end{array}
 $$
which have genus one and their $j$-invariants are $-\frac{117649}{89}$,  $\frac{389017}{16}$, $-\frac{4096}{43}$, $\frac{2146689}{145}$, $\frac{3375}{53}$, $-\frac{4096}{43}$ and $\frac{1771561}{612}$. They
  correspond to the elliptic curves in the statement.
 By the splitting of the jacobians and the  equations of these curves, we  obtain  that  all their automorphism groups have order $2$. It is easy to check  that the automorphism groups of the remaining curves are trivial.
\hfill $\Box$

\begin{prop}\label{genus4}\label{prop11}
Among the curves of genus four $X_0^*(202 )$, $X_0^*(262 )$, $X_0^*(
354)$, $X_0^*(366)$, $X_0^*( 370)$, $X_0^*(426 )$, $X_0^*( 546)$ and
$X_0^*(570 )$,  only $X_0^*(370 )$ is bielliptic. The corresponding
quotient curve  is the elliptic curve labeled as $E370a1$ {{ and the
automorphism group of $X_0^*(370)$ has order 2}}. The automorphism
groups of the curves $X_0^*(202 )$, $X_0^*(262 )$, $X_0^*( 354)$,
$X_0^*(426 )$, $X_0^*( 546)$ and $X_0^*(570 )$  are trivial. The
automorphism group of $X_0^*(366)$ has order 2 and the quotient
curve has genus
 2.
\end{prop}
\noindent{\bf Proof.} The splitting of  $J_0^*(N)$ is:
$$\begin{array}{rrcrrr}
J_0^*(202)\stackrel{\Q}\sim & \prod_{i=1}^2 A_{f_i}  \,, & &A_{f_1}\stackrel{\Q}\sim E101a\,,\, &f_2\in\New_{202}^*\,,\,&\dim A_{f_2}=3\,,\\
J_0^*(262)\stackrel{\Q}\sim &\prod_{i=1}^3 A_{f_i} \,, &A_{f_1}\stackrel{\Q}\sim E131a\,,&\,A_{f_2}\stackrel{\Q}\sim E262b\,,\,\, &f_3\in\New_{262}^*\,,\, &\dim A_{f_3}=2\,,\\
J_0^*(354)\stackrel{\Q}\sim & \prod_{i=1}^3 A_{f_i} \,, & A_{f_1}\stackrel{\Q}\sim E118a\,,\, &A_{f_2}\stackrel{\Q}\sim E354b\,,\,&f_3\in\New_{177}^*\,,\, &\dim A_{f_3}=2\,,\\
J_0^*(366)\stackrel{\Q}\sim &\prod_{i=1}^3 A_{f_i} \,, &A_{f_1}\stackrel{\Q}\sim E61a\,,&\,A_{f_2}\stackrel{\Q}\sim E122a\,,\,\, &f_3\in\New_{183}^*\,,\, &\dim A_{f_3}=2\,,\\
J_0^*(370)\stackrel{\Q}\sim & \prod_{i=1}^4 A_{f_i} \,, & A_{f_1}\stackrel{\Q}\sim E37a\,,\,& A_{f_2}\stackrel{\Q}\sim E185a\,,& A_{f_3}\stackrel{\Q}\sim E185c\,,\,& A_{f_4}\stackrel{\Q}\sim E370a\,,\\
J_0^*(426)\stackrel{\Q}\sim &\prod_{i=1}^3 A_{f_i} \,, &A_{f_1}\stackrel{\Q}\sim E142a\,,&\,A_{f_2}\stackrel{\Q}\sim E426b\,,\,\, &f_3\in\New_{243}^*\,,\, &\dim A_{f_3}=2\,,\\
J_0^*(546)\stackrel{\Q}\sim &  \prod_{i=1}^2 A_{f_i}  \,, & &A_{f_1}\stackrel{\Q}\sim E91a\,,\, &f_2\in\New_{273}^*\,,\,&\dim A_{f_2}=3\,,\\
J_0^*(570)\stackrel{\Q}\sim &\prod_{i=1}^4 A_{f_i} \,, &A_{f_1}\stackrel{\Q}\sim E57a\,,\,&A_{f_2}\stackrel{\Q}\sim E190b\,,\,&A_{f_3}\stackrel{\Q}\sim E285b\,,\, &A_{f_4}\stackrel{\Q}\sim E570a\,.
\end{array}
$$
In all cases, $\dim \cL_2=1$. Next, we show a generator $Q_2(x,y,z,t)\in\cL_2$:
$$\begin{array}{c|r}
N & Q_2 \phantom{ccccccccccccccccccccccccccccccccccccccc}\\ \hline
202 & -9 t^2 + x^2 + 9 t y - 2 x y + y^2 + 9 t z - x z - 8 y z - 2 z^2\\
262&7 t^2 - t x + x^2 - 4 t y - x y - 4 t z - x z + y z\\
354 &-42 t^2 - 5 t x + 2 t y - 10 x y + 22 y^2 + 48 t z - 15 x z + 6 y z -
 3 z^2
\\
 366 & -4 t^2 + x^2 + 2 x y - y^2 - 2 z^2 \\
 370&-144 t^2 + 27 x^2 - 72 x y - 32 y^2 + 54 x z + 40 y z + 127 z^2\\
 426 &-108 t^2 + 165 t x + 100 x^2 + 78 t y + 55 x y - 137 y^2 + 27 t z -
 30 x z + 84 y z - 72 z^2\\
 546 &388 t^2 + 36 t x + x^2 - 68 t y - 10 x y + 9 y^2 - 420 t z - 34 x z +
 58 y z + 125 z^2\\
 570 &-327 t^2 - 232 t x + 176 x^2 + 36 t y + 64 x y + 192 y^2 + 70 t z -
 40 x z + 188 y z - 127 z^2
 \end{array}
$$
Only $X_0^*(370)$ could be bielliptic. In this case, the curve is trigonal (see \cite[Proposition 1]{HaSh}) and $\dim \cL_3=5$.
By computing a polynomial  $Q_3\in\cL_3 $ that is not multiple of $Q_2$, we get
$$ \begin{array}{cr}Q_3(x,y,z,t)=&27 x^3 - 90 x^2 y + 32 y^3 + 63 x^2 z + 108 x y z + 81 x z^2 -
 114 y z^2 - 107 z^3\,.\\
\end{array}
$$
Since $Q_3(x,y,z,-t)\in \cL_3$,  the curve $X_0^*(370)$ is bielliptic by Lemma \ref{con}. Let us check this result.
 Set $P(X,Y):=\operatorname{Resultant\,}(Q_2(X,Y,1,T),Q_3(X,Y,1,T),X)$.
      More precisely,
      $$  \begin{array}{cr}\label{equation}
   P(T,Y)=& -592 + 1944 T^2 - 4860 T^4 + 2916 T^6 - 4752 T^2 Y + 3402 T^4 Y -
 408 Y^2 \\& + 729 T^4 Y^2 + 80 Y^3 + 1620 T^2 Y^3 + 396 Y^4 -
 270 T^2 Y^4 - 222 Y^5 + 17 Y^6 \,.\end{array}
$$
The curve determined by the    equation $ P(T,Y)=0$ has genus $4$.
Hence, it is a plane model for  $X_0^*(370)$. The model admits  the
involution $u\colon (T,Y)\mapsto (-T,Y)$. Replacing $T^2$ with $T$,
we obtain   a genus one curve, whose $j$-invariant is $
15438249/2960$. Checking \cite[Table1]{Cre}, the elliptic quotient
has conductor $370$ and  label $a1$. The polynomials $Q_2$ and $Q_3$
show that  $u$ is the only nontrivial involution of the curve. It is
clear that the remaining curves, except $X_0^*(366)$, have trivial
automorphism group.

Looking at  the polynomial $Q_2$ for $N=366$, we may ask  whether  one of the two  linear maps $(\omega_1,\omega_2,\omega_3,\omega_4) \mapsto \pm (-\omega_1,-\omega_2,\omega_3,\omega_4)$  comes from an involution $u$ of $X_0^*(366)$. After determining $\cL_3$, the answer is affirmative. Hence, $\Jac (X_0^*(366))$ is isogenous over $\Q$ to $E61a\times E122a$ or $A_{f_3}$. After  checking  which of the vector subspaces $\langle \omega_1,\omega_2\rangle$ or $\langle \omega_3,\omega_4\rangle$ provides a hyperelliptic curve,  the right  answer is $A_{f_3}$, and an equation for the quotient curve $X_0^*(366)/\langle u\rangle$ is
$$
Y^2 =X^6 - 6 X^5 + 23 X^4 - 42 X^3 + 53 X^2 - 24 X + 4\,.\qquad \Box
$$

   \begin{rem}\label{rem5}
It is expected that the automorphism group of $X_0^*(N)$ is trivial
for a large enough $N$ and, thus, the genera of the quotients curves
by nontrivial involutions are bounded. The curve $X_0^*(366)$
  shows that if  this bound exists,  then it {{is}} at least $2$.
  \end{rem}

\begin{prop} The curves $X_0^*(402)$, $X_0^*(438)$, $X_0^*(714)$,  $X_0^*(798)$, $X_0^*(910)$, $X_0^*(690)$, $X_0^*(858)$ and  $X_0^*(870)$ are not bielliptic.

\end{prop}
\noindent{\bf Proof.}
The splitting of $J_0^*(N)$ for the  curves of genus $5$ in the statement is:
$$\begin{array}{rrrrrr}
J_0^*(402)\stackrel{\Q}\sim& \prod_{i=1}^4A_{f_i}\,, & A_{f_1}\stackrel{\Q}\sim E201a\,,\,& A_{f_2}\stackrel{\Q}\sim E201c\,,\,&A_{f_3}\stackrel{\Q}\sim E402a\,,\,&f_4\in\New_{67}^*\,,\,\\ &&&&&\dim A_{f_4}=2\,,\\
J_0^*(438)\stackrel{\Q}\sim& \prod_{i=1}^4A_{f_i}\,, & A_{f_1}\stackrel{\Q}\sim E219a\,,\,& A_{f_2}\stackrel{\Q}\sim E219c\,,\,&A_{f_3}\stackrel{\Q}\sim E438a\,,\,&f_4\in\New_{73}^*\,,\,\\ &&&&&\dim A_{f_4}=2\,,\\
J_0^*(714)\stackrel{\Q}\sim& \prod_{i=1}^4A_{f_i}\,, &  A_{f_1}\stackrel{\Q}\sim E102a\,,&A_{f_2}\stackrel{\Q}\sim E238b\,\,,\,&A_{f_3}\stackrel{\Q}\sim E714a \,,\,&f_4\in \New_{357}^*\\ &&&&& \dim A_{f_4}=2\,,\\
J_0^*(798)\stackrel{\Q}\sim &\prod_{i=1}^4A_{f_i}\,, &A_{f_1}\stackrel{\Q}\sim E57a\,,\, &A_{f_2}\stackrel{\Q}\sim E399a\,,\, & A_{f_3}\stackrel{\Q}\sim E798a\,,\,&f_4\in\New_{133}^*\,,\\ &&& &&\dim A_{f_4}=2\,,\\
J_0^*(910)\stackrel{\Q}\sim &\prod_{i=1}^4A_{f_i}\,,
&A_{f_1}\stackrel{\Q}\sim E65a\,,\,&A_{f_2}\stackrel{\Q}\sim
E91a\,,&\,A_{f_3}\stackrel{\Q}\sim
E455a\,,\,&f_4\in\New_{910}^*\,,\,\\& && &
  &\dim A_{f_4}=2\,.
\end{array}
$$
{{ In all cases to study $(N,E)$, we have $\dim \cL_{2,i}=0$, with
$i$ the one corresponding to $E$. More explicitly, for
$N=402,714,798, 910$ we have $\dim \cL_{2,2}=0$, for $N=438$, $\dim
\cL_{2,1}=0$, and for $N=910$,
 also $\dim \cL_{2,3}=0$.}}
\vskip 0.2 cm

For the curves of genus 6, the splitting of  $J_0^*(N)$ is:
$$\begin{array}{rrrrrr}
J_0^*(610)\stackrel{\Q}\sim& \prod_{i=1}^4A_{f_i}\,, & A_{f_1}\stackrel{\Q}\sim E138a\,,\,& A_{f_2}\stackrel{\Q}\sim E690a\,,\,&f_3\in\New_{345}^*\,,\,&f_4\in\New_{115}^*\,,\,\\ &&&&\dim A_{f3}=2&\dim A_{f_4}=2\,,\\
J_0^*(858)\stackrel{\Q}\sim& \prod_{i=1}^4A_{f_i}\,, &  A_{f_1}\stackrel{\Q}\sim E143a\,,&A_{f_2}\stackrel{\Q}\sim E290a\,\,,\,&f_3\in \New_{429}^* \,,\,&f_4\in \New_{858}^*\\ &&&&\dim A_{f_3}=2& \dim A_{f_4}=2\,,\\
\end{array}
$$
In all cases, $\dim \cL_{2,1}=1$. Hence, $\dim \cL_{2,1}=1<\dim \cL_2=6$. For $N=898$,  also $\dim \cL_{2,1}=1$.
\vskip 0.2 cm

Finally,  the splitting  for the curve  $X_0^*(870)$ of genus seven is:
$$\begin{array}{rrrrrr}
J_0^*(870)\stackrel{\Q}\sim&\prod_{i=1}^4A_{f_i}\,, &  A_{f_1}\stackrel{\Q}\sim E58a\,,\,&A_{f_2}\stackrel{\Q}\sim E145a\,,\,& A_{f_3}\stackrel{\Q}\sim E290a\,,\,&f_4\in\New_{435}^*\,,\,\\ &&&&&\dim A_{f_4}=4\,.\\
\end{array}
$$
In this case, $\dim \cL_{2,2}=\dim \cL_{2,3}=4 <\dim \cL_2=10$.
\hfill $\Box$

\vskip 0.2 cm
As a consequence of the previous results, we obtain the statement of Theorem \ref{main} for $N$ even.

\begin{cor} For $N$ even, the curve $X_0^*(N)$ is bielliptic exactly for the thirteen values of $N$ in the set
$$
\{106,122,158,166,178, 246,{{258}},290,318,370,390,430,510\}\,.
$$
For these values of $N$  automorphism group has order
$2$ when $g_N^*>2$, otherwise it is the Klein group.
\end{cor}

\section{Quadratic points}

Let us now  prove Theorem 2. We know by \cite{HaHa} that if $N$ is
square-free and  $X^*_0(N)$ is hyperelliptic, then $g_N^*=2$. On the
other hand, a genus two curve defined over a number field $K$
 is hyperelliptic over $K$ and, thus,
 all genus two curves  $X_0^*(N)$ are hyperelliptic over
$\Q$. The set of values of $N$ in Theorem 2 are those for which
$g_N^*=2$ and those such that  $X_0^*(N)$ is bielliptic and
$g_N^*\geq 3$. This is due to the fact that, when $g_N^*\geq 3$,
the quotient curve is always an elliptic curve with rank equal  to 1
(see \cite[Table1]{Cre}).  Hence, all these values of $N$ are
exactly the values for which  $\Gamma_2( X_0^*(N),\Q)$ is infinite
(cf. \cite[Theorem 2.14]{BaMom}).

\section{Appendix}

Here we list the values $N$ such that
$g_N^*\leq 2$. The table for genus 2 reproduces the one in
\cite{Hata}.  The tables for genus 0 or 1 are taken from
\cite{GL}.   We note that  the value
 $N=141$, which  does not appear in Proposition 1.1 of \cite{GL},  is included in  the appendix of this paper and here.

\vskip 0.2 cm

\begin{tabular}{|c|r|}
\hline $g_N^*=0$&
2,3,5,6,7,10,11,13,14,15,17,19,21,22,23,26,29,30,31,33,34,35,{{38,39,}}\\
&41,42,46,47,51,55,59,62,66,69,70,71,78,87,94,95,105,110,119\,.\\
\hline $g_N^*=1$&
37,43,53,57,58,61,65,74,77,79,82,83,86,89,91,101,102,111,114,118,123,\\
&130,131,138,141,142,143,145,155,159,174,182,190,195,210,222,231,238\,.\\

\hline $g_N^*=2$& {{67}}, 73,85,93,103,106,107,115,122,
129,133,134,146,154,158,161,165, {{166}},167\\ &170,
177,186,191,205,206,209,213,215,221,230,255,266,285,286,287, 299\\
&330,357,390\,.\\ \hline
\end{tabular}

\vskip 0.4 cm {{ \noindent {\bf Acknowledgements.} { We thank the
referees for their comments, especially those that have contributed
to improve the computations of the different tables and the
exposition of the paper. }}}

\bibliographystyle{abbrv}

\noindent{Francesc Bars Cortina}\\
{Departament Matem\`atiques, Edif. C, Universitat Aut\`onoma de Barcelona\\
08193 Bellaterra, Catalonia}\\
{francesc@mat.uab.cat}

\vspace{1cm}

\noindent
{Josep Gonz\'alez Rovira}\\
{Departament de Matem\`atiques, Universitat Polit\`ecnica de Catalunya EPSEVG,\\
Avinguda V\'{\i}ctor Balaguer 1, 08800 Vilanova i la Geltr\'u,
Catalonia}\\
{josep.gonzalez@upc.edu}

\end{document}